\documentclass[11pt]{amsart}

\usepackage[colorlinks=true, pdfstartview=FitV, linkcolor=blue,citecolor=green, urlcolor=blue]{hyperref}
\usepackage{graphicx}
\usepackage{floatflt}

\newcommand{\nind}{{\noindent}}

\newcommand{\vit}{\nu} 								
\newcommand{\pha}{\mathcal{ P } \left( \mathbb{ R }^2 \right)} 	
\newcommand{\vor}{\Upsilon} 							
\newcommand{\Bd}{{\rm Bd}}							
\newcommand{\Int}{{\rm Int}}							
\newcommand{\Sim}{{\rm Sim}}

\newcommand{\R}{{\mathbb R}}

\newcommand{\N}{{\mathbb N}}

\newcommand{\ppp}{{\mathcal P}}
\newcommand{\sss}{{\mathcal S}}

\pdfoutput=1

\newtheorem{dfn}{Definition}
\newtheorem{thm}{Theorem}
\newtheorem{prp}{Proposition}
\newtheorem{cor}{Corollary}
\newtheorem{lem}{Lemma}
\newtheorem*{thm*}{Theorem}

\newcounter{count}

\begin{document}
\title{A Dynamical System Using the Vorono\"{i} Tessellation}
\author{Natalie Priebe Frank \& Sean Hart}
\date{\today}
\address{Department of Mathematics\\Vassar College\\Poughkeepsie, NY  12604\\nafrank@vassar.edu\\sehart@vassar.edu}
\setlength{\baselineskip}{.5cm}

\maketitle

\section{Introduction}
\label{sec:int}
Suppose you know the locations of post offices or cell phone satellites, and you want to know what regions they serve.  Or maybe you know the locations of atoms in a crystal, and you want to know what a fundamental region looks like.  There are lots of reasons you might want to make a tiling around a given discrete set of points.   A natural way to do it is with ``Vorono\"{i} tessellations"---so natural, in fact, that it has been rediscovered numerous times over the years.

On the other hand, if you have a tiling, you might want to decorate each tile with a few points to create or destroy symmetry.  Or you might look at all the {\bf vertices} of the tiling--points where three or more tiles meet--to extract combinatorial information.  
If you have a tiling, there are many ways to extract a point set from it.  

So, you can get tilings from point sets and point sets from tilings:  doesn't this give you a way to associate point sets to point sets or tilings to tilings?  Once you have a map from a class of objects back to itself, you can take a dynamical systems viewpoint to analyze the situation.  In this paper we are going to do exactly that, with a new dynamical system based on the vertices of  Vorono\"{i} tessellations.

For those uninitiated with the Vorono\"{i} tessellation, we begin with its definition and then give the definition of our dynamical system.  From there, the remainder of $ \mathsection $\ref{sec:int} is spent exploring the evolution of simple point sets, using these simplified examples to develop both the intuition and vocabulary needed for more interesting cases.  In $ \mathsection $\ref{sec:whatno} we give a new proof of a theorem, first proved in \cite{efgh}, quantifying the growth in size of point sets over repeated iteration.  Following that we will point out some interesting corollaries and give some estimates on the growth rate.  We
devote $ \mathsection $\ref{sec:dsv} to discussing what questions interest us most from the dynamical systems viewpoint.

For now, let's turn to the definitions.

\subsection{Our dynamical system}
\label{ss:ourds}

\begin{figure}[ht]
\centering
$ \begin{array}{cc}
\includegraphics[width=70mm]{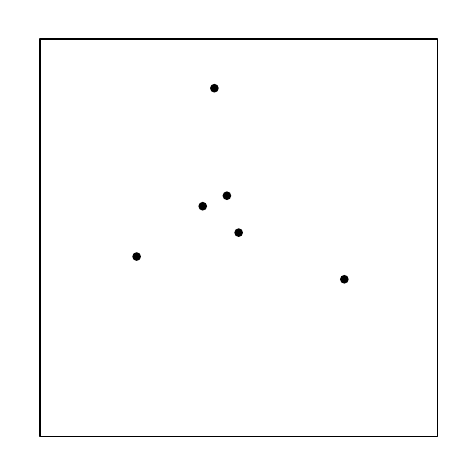} &
\includegraphics[width=70mm]{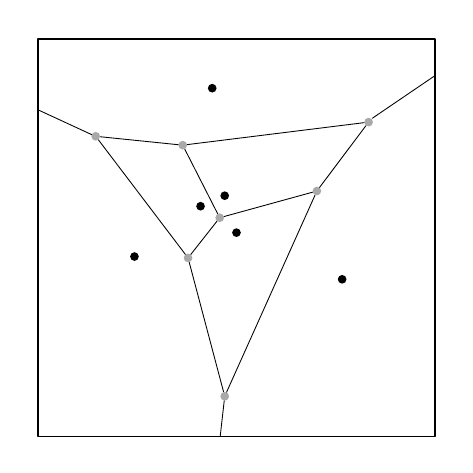} \\
(a) & (b) \\
\end{array} $
\caption{A point set $ P $, it's Vorono\"{i} tessellation, and $ \vit ( P ) $.}
\label{fig:ex1}
\end{figure}

We start with a discrete point set $ P \subset \mathbb{ R }^2 $, which we call the {\bf generating set}, the members of which we refer to as the {\bf generators}.  To avoid technical complications we typically assume that $P$ is finite.   The {\bf Vorono\"{i} polygon} of a point $ p \in P $, denoted $ V ( p ) $, is given by
\begin{equation*}
V ( p ) = \left. \left\{ x \in \mathbb{ R }^2 \right| || p - x || \leq || p' - x || \mbox{ for all } p' \in P \right\} .
\end{equation*}
Simply put, the Vorono\"{i} polygon of $ p $ contains every point in the plane that is closer to $ p $ than to any other member of $ P $, or is equidistant between $ p $ and a nearby generator point.  The {\bf Vorono\"{i} tessellation} of $ P $ is given by
\begin{equation*}
\vor ( P ) = \left\{ V ( p ) \left| p \in P \right\} \right. .
\end{equation*}
A point set and its Vorono\"{\i} tessellation are given in Figure \ref{fig:ex1}.

Given distinct $ p $, $ p' \in P $, the sets $ V ( p ) $ and $ V ( p' ) $ are not necessarily disjoint; in fact, the boundary of each Vorono\"{i} polygon is shared with other Vorono\"{i} polygons.  If $ V ( p ) \cap V ( p' ) $ is a line, ray, or line segment, we call it a {\bf Vorono\"{i} edge} and denote it $ e_{ p , p' } $.   If the intersection of three or more Vorono\"{\i} tiles is a point, we call that point a {\bf Vorono\"{i} vertex}.  The set of all Vorono\"{i} edges and vertices are given by $ \mathcal{ E } ( \vor ( P ) ) $ and $ \mathcal{ V } ( \vor ( P ) ) $, respectively.

 It is useful to notice that $ e_{ p , p' } $ always lies on the perpendicular bisector of the line between $ p $ and $ p' $.   This gives a method for constructing Vorono\"{i} diagrams:  for $ p \in P $, sketch each perpendicular bisector between $ p $ and and another member of $ P $. Then the Vorono\"{i} polygon $ V ( p ) $ is the intersection of all half-planes created by the perpendicular bisectors.  It is also useful to notice that a Vorono\"{\i} vertex is equidistant from the generator points of the tiles it is in.    For proofs of these properties and a wealth of other information about the Vorono\"{i} tessellations, \cite{bible} is an excellent source.

Now, the set $ \mathcal{ V } ( \vor ( P ) ) $ constitutes a point set in its own right, and so one might naturally wonder what its Vorono\"{i} tessellation looks like.  And we need not stop there--the Vorono\"{i} tessellation of $ \mathcal{ V } ( \vor ( P ) ) $ will have a vertex set, too, so how does its Vorono\"{i} tessellation behave?  We have the makings of a dynamical system on the set $ \pha $ of all finite point sets in the plane.

\addtocounter{count}{1}
\begin{dfn} Let $P \in \pha$.
We define the {\bf Vorono\"{i} iteration} of P to be\begin{equation*}
\vit ( P ) = \mathcal{ V } ( \vor ( P ) ) .
\end{equation*}
For $ n = 2, 3, ... $ we define the {\bf n-th Vorono\"{i} iteration} of $ P $ recursively:
\begin{equation*}
\vit^n ( P ) = \mathcal{ V } \left( \vor \left( \vit^{n-1} ( P ) \right) \right) = \vit \left( \vit^{n-1} ( P ) \right) .
\end{equation*}
\end{dfn}

\nind
Figure \ref{fig:ex1} (b) depicts the Vorono\"{i} tessellation of the six points depicted in \ref{fig:ex1} (a), and so the lighter seven points---the vertices of the tessellation---constitute the Vorono\"{i} iteration of the original set.  
Let's begin looking at the simplest cases this dynamical system has to offer.

\subsection{The really trivial cases}
\label{ss:trtc}

We take the convention of using $ | P | $ to denote the {\bf cardinality} of $ P $.  If $ | P | = 1 $, then the single Vorono\"{i} polygon is all of $ \mathbb{ R }^2 $.  Since $ \mathbb{ R }^2 $ is just one big tile, it has no vertices, and so $ \vit^n ( P ) = \O $ for all $ n \geq 1 $.  Then up the ante:  suppose $ P = \{ p_1 , p_2 \} $.  The Vorono\"{i} tessellation then fractures the plane into two half-planes, split along the perpendicular bisector of $ \overline{ p_1 p_2 } $, the line segment joining $ p_1 $ to $ p_2 $.  The vertex set, however, is still empty, and so $ \vit^n ( P ) = \O $ for all $ n \geq 1 $.

There are two possibilities when $ P $ has three points:  either all three are collinear, or they lie on a circle.  The former yields $ \vit^n ( P ) = \O $ for all $ n \geq 1 $ and in the latter, $ \vit ( P ) $ is a single point.  No matter how large $ P $ is, the special cases of collinearity and cocircularity always work out like this:

\addtocounter{count}{1}
\begin{prp}
All the points in $ P $ are collinear iff $ \vit^n ( P ) = \O $ for all $ n \geq 1 $.  All the points in $ P $ are cocircular iff $ | \vit ( P ) | = 1 $ and $ \vit^n ( P ) = \O $ for $ n > 1 $.
\label{prp:collincocir}
\end{prp}

\nind
The first fact follows naturally from the observation that if $ p $, $ q \in P $ such that $ e_{ p , q } \in \mathcal{ E } \left( \vor ( P ) \right) $, then $ e_{ p , q } $ lies on the perpendicular bisector of the line segment $ \overline{ p q } $.  When all the points are collinear, the perpendicular bisectors are all parallel and thus do not intersect to produce new vertices.  Conversely, if $ \vit ( P ) = \O $ then no Vorono\"{i} edges intersect, which implies all Vorono\"{i} edges are parallel, and hence $ P $ must be collinear.  

To prove the cocircularity result, we introduce the concept of an {\bf empty circle}.   This is a circle whose interior does not contain any generator points and whose boundary passes through three or more generator points.  Such a circle gets its name since it is ``empty'' of generator points.   Proposition \ref{prp:collincocir} follows immediately from the next Proposition.

\addtocounter{count}{1}
\begin{prp}
A point $ q \in \mathbb{ R }^2 $ is a Voronoi vertex in $ \mathcal{ V } \left( \vor ( P ) \right) $ if and only if it is the center of an empty circle.\footnote{See \cite{bible}, pg. 61.}
\label{prp:vertcirc}
\end{prp}


\nind
In addition to proving Proposition \ref{prp:collincocir}, this proposition also lets us quickly find the vertex set of $\vor( P) $.  To do so, pick a noncollinear triple $ p $, $ q $, $ r \in P $, and look at the unique circle passing through them.  If this circle is empty, then place a vertex at its center.  Once you have checked every triple, every vertex will be accounted for---and you never had to sketch an edge!  

So, we know what happens when the point set is really small, or collinear, or cocircular.   Let's move on to...

\subsection{A slightly less trivial case}
\label{ss:sltc}

Suppose $ | P | = 4 $.  To discriminate between configurations that aren't collinear or cocircular, we require a more sophisticated vocabulary.  We say that a subset $ A \subset \mathbb{ R }^2 $ is {\bf convex} if, given any distinct $ x $, $ y \in A $, the line segment $ \overline{ x y } $ is contained in $ A $.  
(For practice, try proving that Vorono\"{\i} polygons are convex!)
The {\bf convex hull} of a set $ P $ is defined to be the smallest convex set containing $ P $.  We write $ CH ( P ) $ to mean the convex hull of $ P $, and write $ \partial CH ( P ) $ to restrict attention to its boundary.
\addtocounter{count}{1}
\begin{dfn}
For $ p \in P $, we say $ p $ is on the {\bf boundary of} $ P $ and write $ p \in \Bd ( P ) $ if $ p \in P \cap \partial CH ( P ) $.  Otherwise, we say $ p $ is in the {\bf interior of} $ P $ and write $ p \in \Int ( P ) $.
\label{dfn:bdint}
\end{dfn}
\nind
In general it is clear that
\addtocounter{count}{1}
\begin{equation}
| \Int ( P ) | + | \Bd ( P ) | = | P | .
\end{equation}

We say that two distinct points $ p, p' \in \Bd ( P ) $ are {\bf neighbors on the boundary of} $ P $ if $ \overline{ p p' } \subset \partial CH ( P ) $ and no $ p'' \in P - \{ p , p' \} $ lies on $ \overline{ p p' } $.  For noncollinear finite point sets, any point on the boundary has two distinct neighbors.  

It is not surprising that the Vorono\"{\i} tiles of boundary points would have infinite edges.   We denote by
$\mathcal{ E }_F ( \vor ( P ) ) $  and $\mathcal{ E }_I ( \vor ( P ) ) $ the finite and infinite edges of $\vor(P)$, respectively.
 

\addtocounter{count}{1}
\begin{prp}
Assume not all $ p \in P $ are collinear, and let $ e_{p, p'} \in \mathcal{ E } \left( \vor ( P ) \right) $.  Then $ e_{ p , p' } \in \mathcal{ E }_I \left( \vor ( P ) \right) $ if and only if $ p $ and $ p' $ are neighbors on the boundary of $ P $.  Moreover, $ e_{p, p'} \in \mathcal{ E }_F \left( \vor ( P ) \right) $ if and only if $ p $ and $ p' $ are not neighbors on the boundary of $ P $.\footnote{Follows from \cite{bible}, pg. 59.}
\label{prp:convedge}
\end{prp}

\nind
This then implies
\addtocounter{count}{1}
\begin{equation}
\left| \mathcal{ E }_I \left( \vor ( P ) \right) \right| = | \Bd ( P ) | .
\end{equation}

\begin{figure}[ht]
\centering
$ \begin{array}{cc}
\includegraphics[width=70mm]{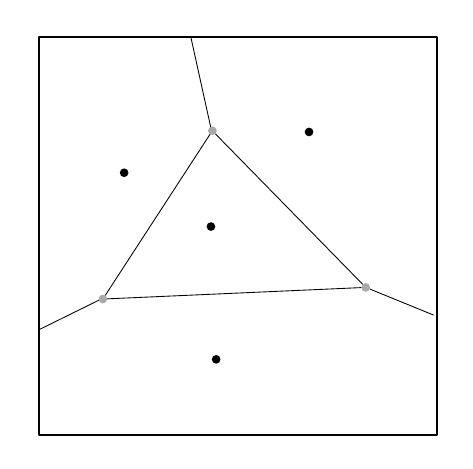}
&
\includegraphics[width=70mm]{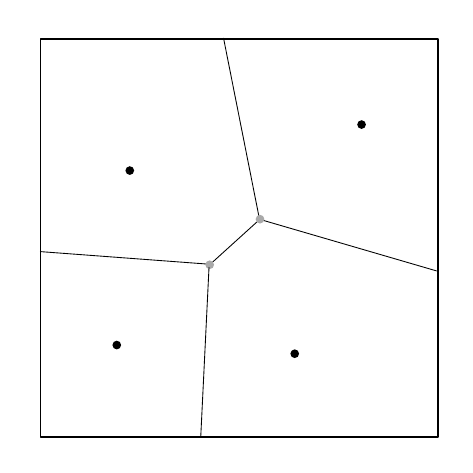} \\
(a) & (b) \\
\end{array} $
\caption{Configurations satisfying (a) $ | \Bd ( P ) | = 3 $ and (b) $ | \Bd ( P ) | = 4 $.}
\label{fig:threeone}
\end{figure}

Returning to the case of $ | P | = 4 $, let us assume that the points in $P$ are neither collinear nor cocircular.   We may have $ | \Bd ( P ) | = 3 $ or $ 4 $, and
we are going to prove that Figure \ref{fig:threeone} represents the only possible types of Vorono\"{\i} iterations.

First let's assume $ | \Bd ( P ) | = 3 $.  Since $ | \Int ( P ) | + | \Bd ( P ) | = | P | = 4 $, say $ \Int ( P ) = \{ p \} $.  Proposition \ref{prp:convedge} tells us that each of the edges of its Vorono\"{i} polygon $ V ( p ) $ must be finite.  As $ | P - \{ p \} | = 3 $, there can be at most three edges; since $ V ( p ) $ must be bounded, there must be exactly three edges.  So $ V ( p ) $ is a triangle, and the three infinite edges promised by Proposition \ref{prp:convedge} extend from its vertices.  Since $|\vit(P)|=3$, we conclude that $ | \vit^2 ( P ) | = 1 $ and $ \vit^n ( P ) = \O $ for all $ n \geq 3 $.

Now let's assume that $|\Bd(P)| = 4$ and that $ P $ is neither collinear nor cocircular.  We will prove that $ | \vit ( P ) | = 2 $.  By Proposition \ref{prp:vertcirc} we know that $ | \vit ( P ) | >1 $.   To prove that $ | \vit ( P ) | \le 2 $, we use a trick that will come in handy again, during the proof of our main theorem.

For $ q \in \vit ( P ) $, let $ \rho ( q ) $ be the {\bf degree} of the vertex $ q $:  the number of edges touching $ q $.  Then we have $ \rho ( q ) \geq 3 $ for each $ q \in \vit ( P ) $.  Summing the degrees of all the vertices therefore is greater than or equal to $ 3 \cdot | \vit ( P ) | $.  Counting edges instead we see that each infinite edge touches exactly one vertex but each finite edge touches exactly two vertices.  Thus
\begin{equation*}
3 \cdot | \vit ( P ) | \leq \sum_{q \in \vit ( P )} \rho ( q ) = \left| \mathcal{ E }_I \left( \vor ( P ) \right) \right| + 2 \cdot \left| \mathcal{ E }_F \left( \vor ( P ) \right) \right| .
\end{equation*}
Since there are only two pairs of points in $ P $ that are not neighbors, Proposition \ref{prp:convedge} implies $ \left| \mathcal{ E }_F \left( \vor ( P ) \right) \right| \leq 2 $.  Thus
\begin{equation*}
3 \cdot | \vit ( P ) | \leq 4 + 2 \cdot 2 = 8 ,
\end{equation*}
and so $ | \vit ( P ) | \leq 2 $, proving that in this case, $ | \vit ( P ) | = 2 $.
Moreover we can conclude that $ \vit^n ( P ) = \O $ for all $ n > 1 $.

\subsection{The rest of the cases are all hard.}
\label{ss:allhard}

\begin{figure}[ht]
\centering
$ \begin{array}{cc}
\includegraphics[width=70mm]{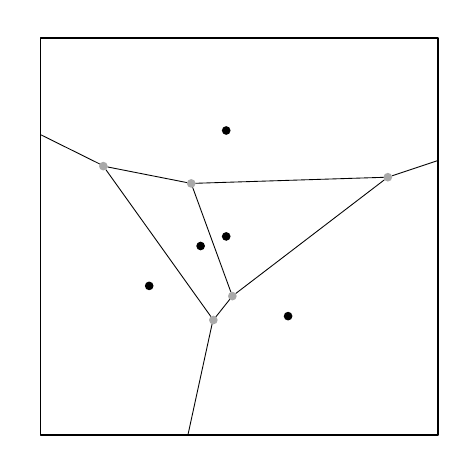} &
\includegraphics[width=70mm]{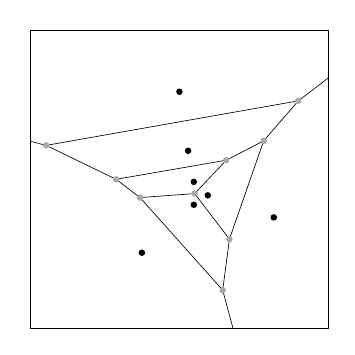} \\
(a) & (b) \\
\end{array} $
\caption{Configurations for which $ | \vit ( P ) | $ is (a) equal to or (b) greater than $ | P | $.}
\label{fig:fivetofive}
\end{figure}

The case $ | P | = 4 $ turns out to be the last case for which $ \vit^n ( P ) $ is guaranteed to equal the null set for large enough $ n $.  Indeed, the case $ | P | = 5 $ has resisted our attempts to understand it!  Figure \ref{fig:fivetofive} depicts configurations of larger cardinality where the point set does not iterate to a simpler configuration, which in turn enriches (and complicates) matters. We must suspend our case-by-case analysis at this point.

\section{How big is $ \vit ( P  )$?}
\label{sec:whatno}
\setcounter{count}{0}

Having played with the dynamical system for a while, we, along with our colleagues Allison Edgren and Olivia Gillham, decided to explore the question of what happens to the size of $  \vit^n ( P )  $ as $ n $ goes to infinity.   We had some success:  a formula \cite{efgh} (unpublished) that tells us exactly how many points will be in $\vit(P)$!  We are going to share this result with you as soon as we describe one more essential piece of geometric information.

\begin{figure}[ht]
\centering
$ \begin{array}{cc}
\includegraphics[width=70mm]{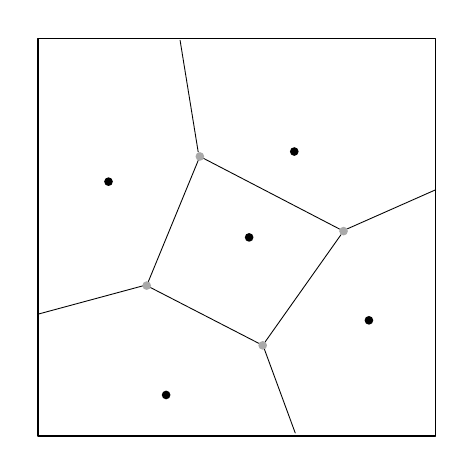} &
\includegraphics[width=70mm]{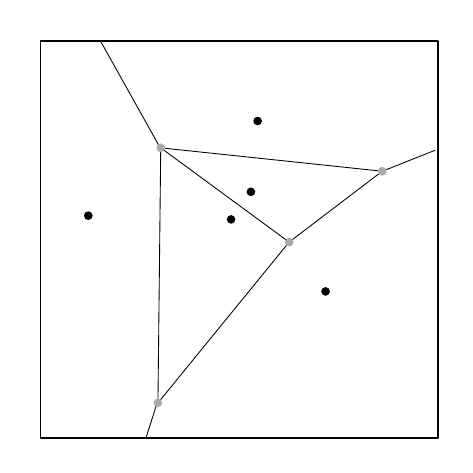} \\
(a) & (b) \\
\end{array} $
\caption{Two configurations of $ 5 $ points.}
\label{fig:illuex}
\end{figure}

\subsection{Relating degeneracy \& cocircularity}
\label{ss:reldc}
Compare the configuration given in Figure \ref{fig:fivetofive} (a) with those sketched in Figure \ref{fig:illuex}.  Each are configurations of five points, yet each has a distinctly different iteration.  We already have the vocabulary necessary to distinguish the configuration sketched in Figure \ref{fig:illuex} (a) from the others:  they differ with respect to the number of points on the boundary.  The difference between the other two, however, is more subtle.

A meticulous reader might notice that all of the vertices in Figure \ref{fig:fivetofive} (a) are of degree three while there is a vertex in \ref{fig:illuex} (b) with degree four.  A Vorono\"{i} diagram with a vertex of degree greater than three is called {\bf degenerate}; if all the vertices of the Vorono\"{i} diagram are degree three we call it {\bf non-degenerate}.  We can quantify this concept of degeneracy with the number $ I_c ( P ) $:  

\addtocounter{count}{1}
\begin{dfn}
Let $ \left\{ C_1 , \dots , C_k \right\} $ be the set of all empty circles of a point set $ P $.  The number of {\bf instances of cocircularity} of $ P $ is given by
\begin{equation*}
I_c ( P ) = \sum_{ i = 1 }^k \left( | C_i \cap P | - 3 \right) .
\end{equation*}
\label{dfn:mtcirc}
\end{dfn}

\nind
In a non-degenerate point set, every vertex will be of degree three, and so every empty circle will intersect exactly three points in $ P $.  In such a case, we compute $ I_c ( P ) = 0 $.  When a vertex $ q $ is of degree $ k > 3 $, this implies that $ | C_q \cap P | = k $.  In this case, we say that $ q $ contributes $ ( k - 3 ) $ instances of cocircularity.  This argument proves the following:

\addtocounter{count}{1}
\begin{prp}
Let $ \rho ( q ) $ denote the degree of a vertex $ q \in \vit ( P ) $.  Then
\begin{equation*}
I_c ( P ) = \sum_{q \in \vit ( P )} \left( \rho ( q ) - 3 \right) .
\end{equation*}
\label{prp:icp}
\end{prp}

\nind
When given a degenerate configuration, one can obtain a non-degnerate point set by simply shifting any cocircular points by an arbitrarily small amount; further, when given a non-degenerate configuration, it remains non-degenerate under sufficiently small perturbations.  Hence, when studying or applying Vorono\"{i} tessellations it is customary to ignore degenerate configurations:  the modified non-degenerate configuration is usually ``close enough'' to the original, and if the points are collected from real-world instruments we can never assume that our measurements are exact anyway.

\begin{figure}[ht]
\centering
$ \begin{array}{cc}
\includegraphics[width=70mm]{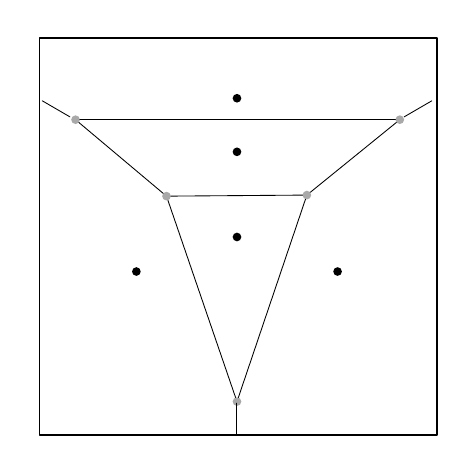} &
\includegraphics[width=70mm]{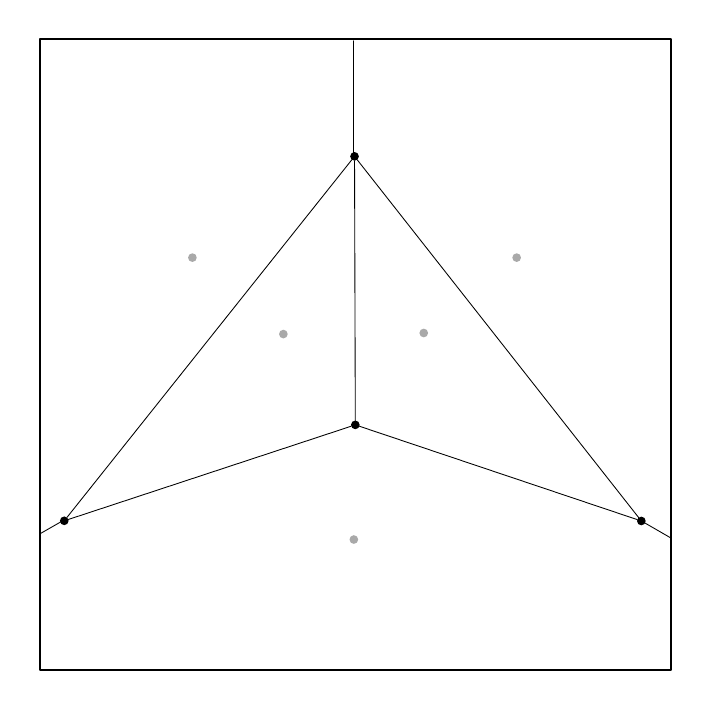} \\
(a) & (b) \\
\end{array} $
\caption{$ I_c ( P ) = 0 $ while $ I_c ( \vit ( P ) ) = 1 $. }
\label{fig:nontodeg}
\end{figure}

If only we could make such an assumption!  In our system, computing $ I_c ( P ) $ is usually only possible if you already know where the vertices are and, as we shall see, incorporating degeneracy makes determining the long-term behavior much more difficult.  But we cannot simply throw out these configurations; there exist situations, like that which is depicted in Figure \ref{fig:nontodeg}, that ruin things for us.  A non-degenerate point set may iterate to a degenerate one, and so assuming $ P $ is non-degenerate does not guarantee $ \vit ( P ) $ is as well.

\nind

\subsection{The theorem on counting vertices.}
\label{ss:countvert}

We are ready to state and prove the theorem discovered by the second author and his colleagues in \cite{efgh} (unpublished).  The proof we present here, due only to the second author, is much simpler than the original.

\addtocounter{count}{1}
\begin{thm}
If not all $ p \in P $ are collinear, then
\begin{equation*}
| \vit ( P ) | = 2 \cdot | P | - | \Bd ( P ) | - I_c ( P ) - 2 .
\end{equation*}
\label{thm:quanti}
\end{thm}

\begin{proof}
We can use the trick from $ \mathsection $\ref{ss:sltc}, summing the degrees of the vertices in $ \vit ( P ) $.  Since $ P $ is non-collinear, every infinite edge intersects exactly one vertex and every finite edge intersects exactly two.  Hence, recalling from $ \mathsection $\ref{ss:sltc} that $ \left| \mathcal{ E }_I \left( \vor ( P ) \right) \right| = | \Bd ( P ) | $, we get
\begin{equation*}
\sum_{ q \in \vit ( P ) } \rho ( q ) = \left| \mathcal{ E }_I \left( \Upsilon ( P ) \right) \right| + 2 \cdot \left| \mathcal{ E }_F \left( \Upsilon ( P ) \right) \right| = | \Bd ( P ) | + 2 \cdot | \mathcal{ E }_F ( \vor ( P ) ) | .
\end{equation*}
Euler's formula for finite planar graphs\footnote{See \cite{henle}, p. 5} tells us that $ V - E + F = 1 $.  We can apply Euler's formula to $ \vor ( P ) $ by ignoring the infinite edges and infinite Vorono\"{i} polygons.  Doing this yields $ | \vit ( P ) | - \left| \mathcal{ E }_F \left( \vor ( P ) \right) \right| + | \Int ( P ) | = 1 $, and so $ \left| \mathcal{ E }_F \left( \Upsilon ( P ) \right) \right| = | \vit ( P ) | + | \Int ( P ) | - 1 $.  Thus
\begin{equation*}
\sum_{ q \in \vit ( P ) } \rho ( q ) = | \Bd ( P ) | + 2 \cdot \left( | \vit ( P ) | + | \Int ( P ) | - 1 \right) .
\end{equation*}
Recalling from $ \mathsection $\ref{ss:sltc} that $ | P | = | \Bd ( P ) | + | \Int ( P ) | $, we have
\begin{equation*}
\sum_{ q \in \vit ( P ) } \rho ( q ) = 2 \cdot | \vit ( P ) | + 2 \cdot | P | - | \Bd ( P ) | - 2 .
\end{equation*}
On the other hand, by Proposition \ref{prp:icp}, 
\begin{equation*}
\sum_{q \in \vit ( P )} \rho ( q ) = 3 \cdot | \vit ( P ) | + I_c ( P ) .
\end{equation*}
Thus
\begin{equation*}
2 \cdot | \vit ( P ) | + 2 \cdot | P | - | \Bd ( P ) | - 2 = 3 \cdot | \vit ( P ) | + I_c ( P )
\end{equation*}
and so, solving for $ | \vit ( P ) | $, we find
\begin{equation*}
| \vit ( P ) | = 2 \cdot | P | - | \Bd ( P ) | - I_c ( P ) - 2 ,
\end{equation*}
as desired.
\end{proof}

\nind
From this and our preliminary analysis in $ \mathsection $\ref{sec:int}, several interesting observations may be made:

\addtocounter{count}{1}
\begin{cor}
\mbox{ }
\begin{enumerate}
\item
If $ | P | < 5 $ then $ \vit^n ( P ) = \O $ for some $ n \in \N$.
\item
If $ | P | = 5 $ then either $ \vit^n ( P ) = \O $ for some $ n \in \N $ or $ | \vit^n ( P ) | = 5 $ for all $ n \in \N $.
\item
If $ | \vit ( P ) | > | P | $ then $ | P | > 5 $.
\item
If $ | \Int ( P ) | = 2 $ and $ I_c ( P ) = 0 $, then $ | \vit ( P ) | = | P | $.
\end{enumerate}
\end{cor}


\nind
Theorem \ref{thm:quanti} alone is not enough to predict the long-term behavior of our system, since computing $ | \vit^n ( P ) | $ requires that we can also compute $ | \Bd ( \vit^{n-1} ( P ) ) | $ and $ I_c ( \vit^{n-1} ( P ) ) $.  We have  been able to find bounds on the size of $\vit^n(P)$, and we present one of interest next.

\subsection{Bounds on the size of $\vit^n (P)$ }
\label{ss:expbd} 

With no assumptions on the geometry of $\vit^n(P)$ we can make an upper bound on the size of $ \vit^{n} ( P ) $, which we conjecture to be sharp.

\addtocounter{count}{1}
\begin{prp}
For all $ n > 0 $, we have $ | \vit^n ( P ) | \leq 2^n \left( | P | - 5 \right) + 5 $.
\end{prp}

\begin{proof}
By induction on $ n $.  For $ n = 1 $, plugging the inequalities $ | \Bd ( P ) | \geq 3 $ and $ I_{c} ( P ) \geq 0 $ into Theorem \ref{thm:quanti} immediately gives us what we want.  Proceeding inductively, we get
\begin{equation*}
| \vit^n ( P ) | = 2 \cdot | \vit^{n-1} ( P ) | - | \Bd ( \vit^{n-1} ( P ) ) | - I_c ( \vit^{n-1} ( P ) ) - 2 \leq 2 \cdot | \vit^{n-1} ( P ) | - 5
\end{equation*}
\begin{equation*}
\leq 2 \cdot \left( 2^{n-1} \cdot \left( | P | - 5 \right) + 5 \right) - 5 = 2^n \cdot \left( | P | - 5 \right) + 5 ,
\end{equation*}
which completes the proof.
\end{proof}

\nind

We have also been able to find an upper bound on the size of the boundary of $\vit(P)$.

\addtocounter{count}{1}
\begin{thm}
The number of points on the boundary of $ \vit ( P ) $ does not exceed $ | \Bd ( P ) | $.
\label{thm:bdvp}
\end{thm}

\nind
We shall need the following lemma:

\addtocounter{count}{1}
\begin{lem}
The interior angle for any vertex in a Vorono\"{i} polygon is strictly less than $ \pi $.
\label{lem:vpqangle}
\end{lem}

\begin{proof}
Let $ p \in P $ generate the Vorono\"{i} polygon $ V ( p ) \in \vor ( P ) $.  Any vertex of $ V ( p ) $ lays at the intersection of some two Vorono\"{i} edges $ e_{p, p'} $ and $ e_{p, p''} $, generated by distinct $ p' $, $ p'' \in P $.  From the argument given in $ \mathsection $\ref{ss:sltc}, we deduce that $ V ( p ) $ is convex, and hence immediately have that the interior angle made by these edges must be less than or equal to $ \pi $.  Seeking a contradiction, assume that this angle is equal to $ \pi $.  By the perpendicular bisector property from $ \mathsection $\ref{ss:trtc}, we may obtain $ p' $ by reflecting $ p $ across $ e_{p, p' } $ and $ p'' $ by reflecting $ p $ across $ e_{p, p''} $.  But then $ p' = p'' $, contradicting the assumption that $ p' $ and $ p'' $ be distinct.
\end{proof}

\begin{proof}[Proof of Theorem \ref{thm:bdvp}]
If all $ q \in \vit ( P ) $ are collinear, then the theorem is trivial.  To prove the result for non-collinear vertex sets, we seek to show $ \left| \Bd \left( \vit ( P ) \right) \right| \leq \left| \mathcal{ E }_I \left( \vor ( P ) \right) \right| $ by proving that every $ q \in \Bd \left( \vit ( P ) \right) $ intersects at least one infinite Vorono\"{i} edge.  This, with the equality $ | \mathcal{ E }_I \left( \vor ( P ) \right) | = | \Bd ( P ) | $ from Proposition \ref{prp:convedge}, gives us what we want.

Since $ q \in \Bd \left( \vit ( P ) \right) $, we may find neighbors $ q' $, $ q'' \in \Bd ( \vit ( P ) ) $ of $ q $ on the convex hull of $ \vit ( P ) $.  Define $ H ( q , q' ) $ to be the closed half-plane through $ q $ and $ q' $ such that $ \vit ( P ) $ is contained in $ H ( q , q' ) $, and similarly define $ H ( q , q'' ) $; these half-planes exist by the properties of the convex hull.  We thus find that $ \vit ( P ) $ is contained in the intersection $ H ( q , q' ) \cap H ( q , q'' ) $, which we denote $ H $.

Let $ V ( p ) \in \vor ( P ) $ be such that $ q \in V ( p ) $ and $ V ( p ) \cap \left( \mathbb{ R }^2 - H \right) \neq \O $.  Let $ e \in \mathcal{ E } \left( \vor ( P ) \right) $ be an edge of $ V ( p ) $ that touches $ q $.  If $ e $ is infinite, then we are done.  If $ e $ terminates in a vertex, then by the convexity of $ H $ we have that $ e $ is contained in $ H $.  Let $ e' \in \mathcal{ E } \left( \vor ( P ) \right) $ be the other edge of $ V ( p ) $ that touches $ q $.  By Lemma \ref{lem:vpqangle}, the interior angle made at $ q $ must be strictly less than $ \pi $.  Since $ V ( p ) $ must have a non-trivial intersection with $ \mathbb{ R }^2 - H $, we find that $ e' $ cannot be contained in $ H $.  Thus it cannot terminate in a vertex, and instead is infinite, as desired.
\end{proof}

\nind
Note that the size of the boundary does not steadily decrease to three in all cases.   There seem to be configurations with boundaries whose sizes stay stable over many iterations.   

Theorems \ref{thm:quanti} and \ref{thm:bdvp} combine to give a lower bound on the size of $\vit^n ( P )$ in the generic situation where there are no instances of cocircularity.

\addtocounter{count}{1}
\begin{thm}
Let $P \in \ppp(\R^2)$ and $N \in \N$.   Suppose that for all $n = 0, 1, 2, ... N-1$, $I_c(\vit^n(P))=0$. Then $ | \vit^N ( P ) | \geq 2^N |P| - (2^N - 1)\left( |\Bd(P)| + 2 \right) $.
\label{thm:lowerbound}
\end{thm}

\begin{figure}
\centering
$ \begin{array}{cc}
\includegraphics[width=72mm]{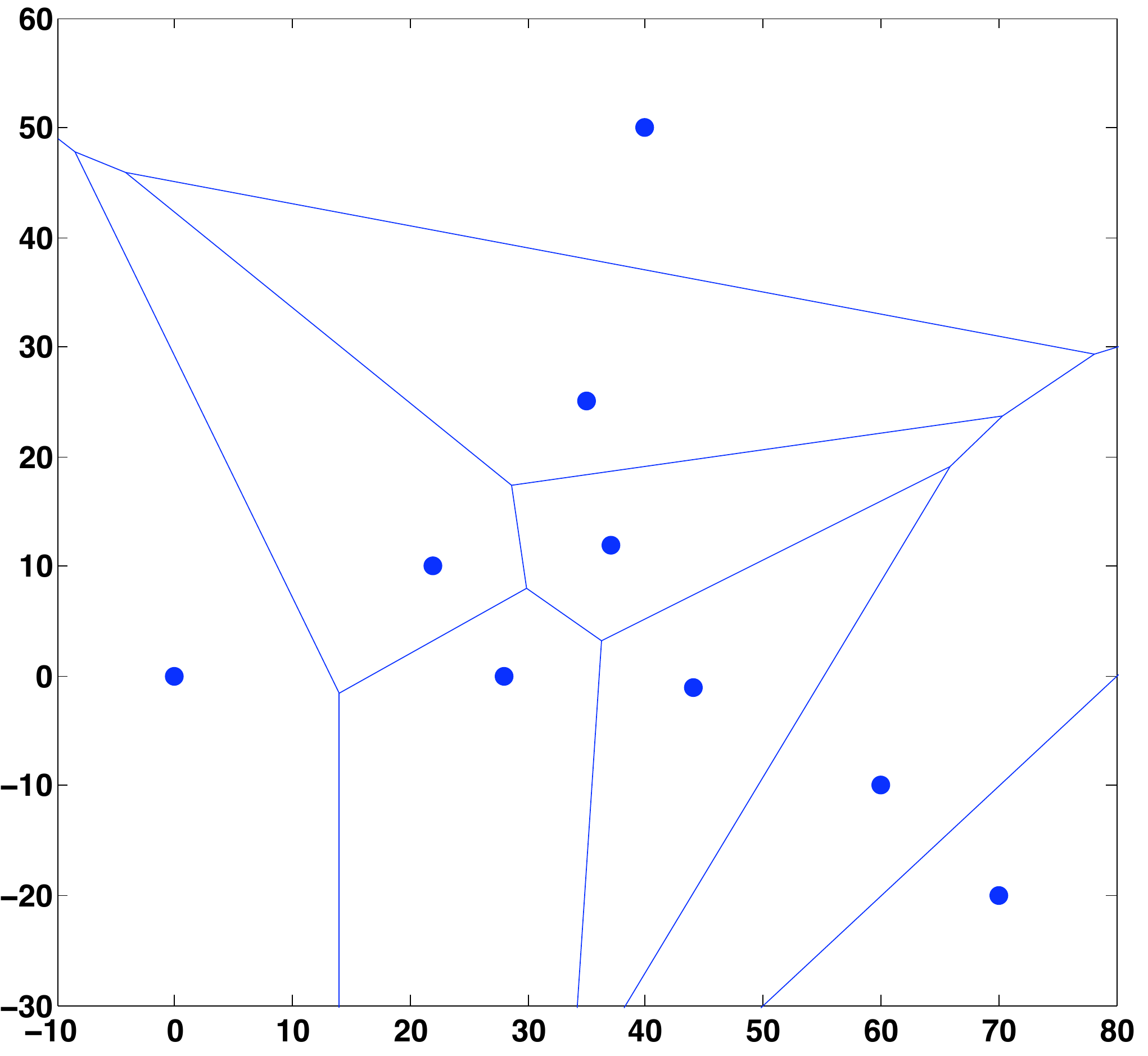} &
\includegraphics[width=70mm]{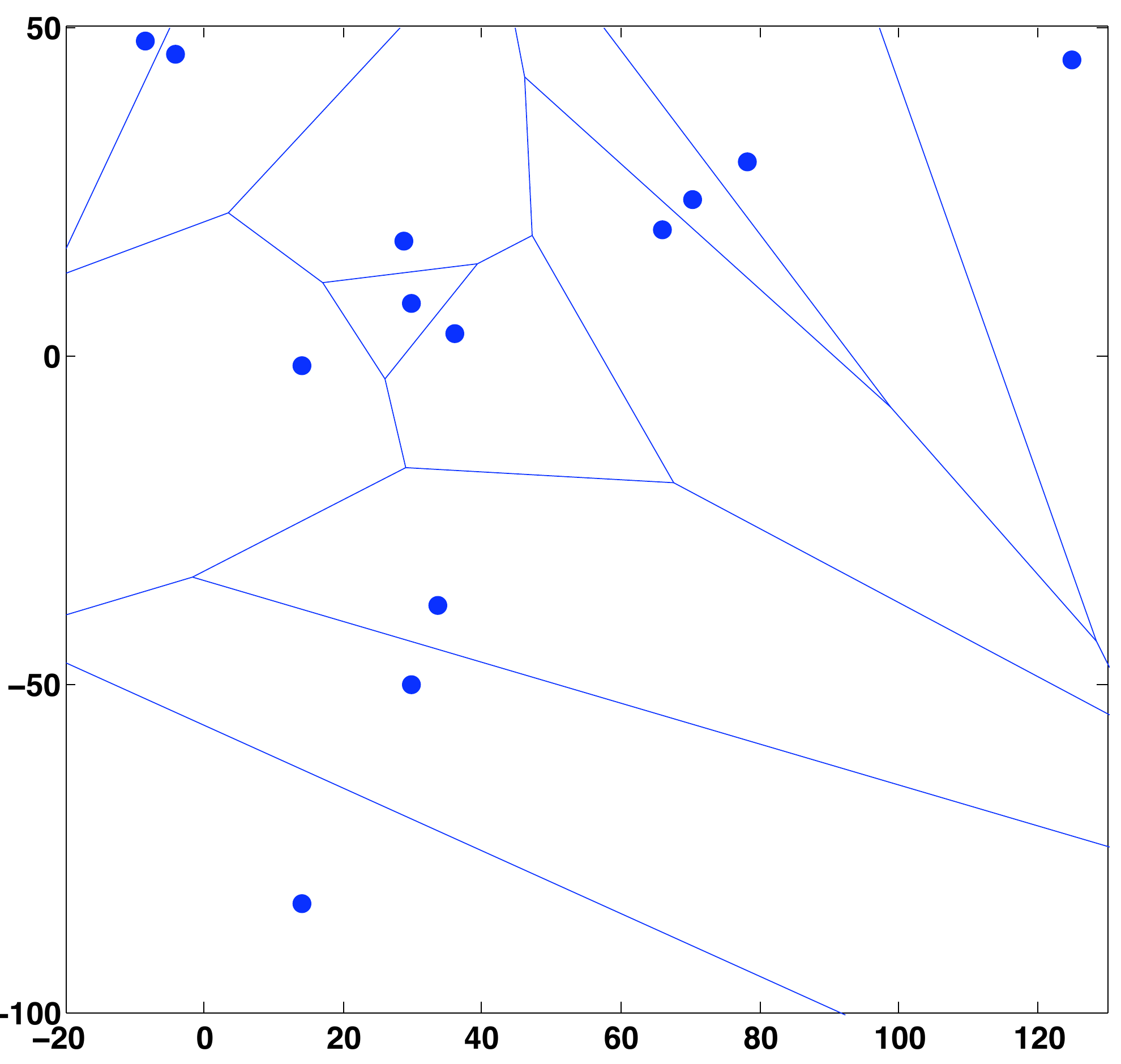} \\ 
\end{array} $
\centerline{\includegraphics[width=80mm]{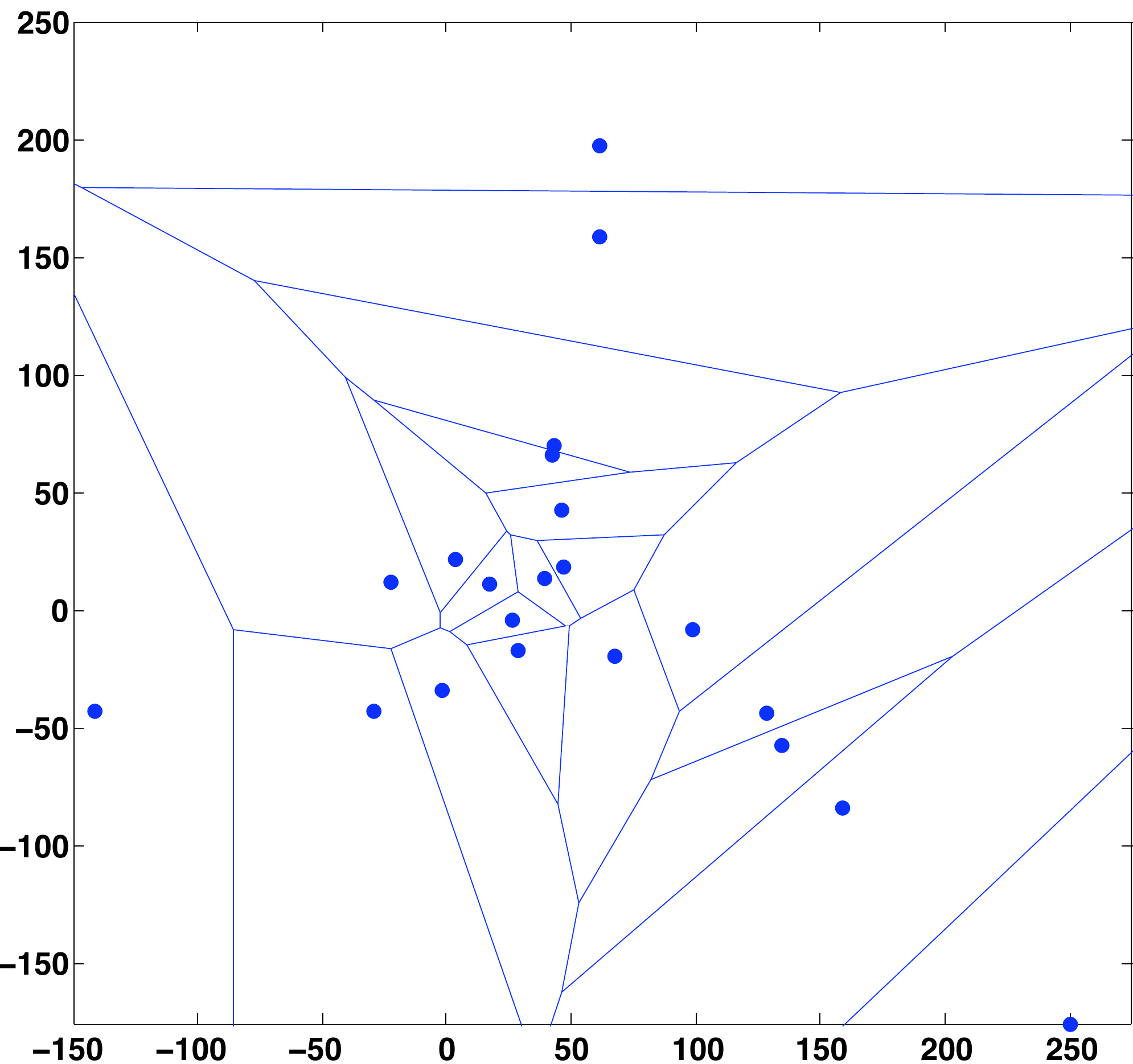} }
\caption{A few iterations of a point set with no instances of cocircularity}
\label{fig:iterationsmessy}
\end{figure}

A point set $P$ with  $|\Int(P)| > 2$ and with Vorono\"{\i} iterations that never contain instances of cocircularity will have exponential growth on the order of $2^n$.   (If $|\Int(P)| \le 2$ the size does not increase).
Since having instances of cocircularity is a rare occurrence, we conjecture that such a point set exists.  We show a typical situation in Figure \ref{fig:iterationsmessy}: a few iterates of a randomly-selected point set of size 9.   As Theorem \ref{thm:quanti} predicts, we see the size going to 13, then 21 points; if there continues not to be any cocircularity the growth will escalate rapidly: the next four iterations contain 37, 69, 133, and then 261 points, respectively.

\section{So what's next? Questions from the dynamical systems viewpoint}
\label{sec:dsv}

In dynamical systems terminology, $ \pha $ is called the ``phase space'' of the system, and elements of $ \pha $ are ``states'' of the system that evolve over time according to the map $ \vit $. 
The {\bf orbit} or {\bf trajectory} of a state $ P $ is defined to be the sequence $ \left\{ P , \vit ( P ) , \vit^2 ( P ) , ... \right\} $ (\cite{chaos} or \cite{katok} are good references).   So far we have only talked about one property of an orbit: the growth rate of $|\vit^n( P )|$ over time.   But there are many other interesting questions we can ask about this dynamical system. 

\subsection{What is the right topology?}
In order to use the machinery of modern dynamics, we ought to have a topology on $\ppp(\R^2)$ making Vorono\"{\i} iteration continuous, at least some of the time.   We  need this to study major dynamical features such as {\em recurrence}---how orbits return to open sets over time---or topological {\em entropy}---a measure of the tendency of the system to become disordered.   Finding the right topology is of great importance, but it has proved to be an interesting problem in its own right. We will mention some of the subtlety here.

To begin, notice that the scale of the diagrams in Figure \ref{fig:iterationsmessy} increases with each iteration.   We want our topology to care about the shape of the Vorono\"{\i} polygons rather than their size.  Similar triangles produce similar Vorono\"{\i} diagrams, so the topology ought to respect that.   Let us be more precise.

\addtocounter{count}{1}
\begin{dfn}
A {\bf similarity transformation} $t \in \Sim(2)$ is a map of the form
\begin{equation*}
t ( {\bf x} ) = k U {\bf x} + {\bf x}_0 ,
\end{equation*}
where $ {\bf x}, {\bf x}_0 \in \mathbb{ R }^2 $, $ U $ is a $ 2 \times 2 $ orthogonal matrix, and $ k \in \mathbb{ R }^+ $. For $ P $, $ Q \subseteq \ppp(\R^2) $, we say $ P $ is {\bf similar} to $ Q $, written $ P \simeq Q $, if there exists $ t \in \Sim ( 2 ) $ such that $ t ( P ) = Q $.
\end{dfn}

\nind
A similarity transformation is composed of translations, rotations, reflections, or dilatations of Euclidean space.  One can show that $ \Sim ( 2 ) $ forms a group under the composition of functions and that the relation $ \simeq $ is an equivalence relation.  Similarity transformations play nice with our dynamical system:

\addtocounter{count}{1}
\begin{thm}
For $ t \in \Sim ( 2 ) $, we have $ t ( \vit ( P ) ) = \vit ( t ( P ) ) $.
\label{thm:simtrans}
\end{thm}
\begin{proof}
Let $ q \in \vit ( P ) $.  By Proposition \ref{prp:vertcirc}, there is a unique empty circle $ C_q $, centered at $ q $.  Since similarity transformations preserve circles, $ t $ maps $ C_q $ to an empty circle centered at $ t ( q ) $.  Thus $ t ( q ) \in \vit \left( t ( P ) \right) $, and so $ t \left( \vit ( P ) \right) \subset \vit \left( t ( P ) \right) $.  By similar logic, since $ t ( P ) $ is a point set and $ t^{-1} $ is a similarity transformation, we have
\begin{equation*}
t^{-1} \left( \vit \left( t ( P ) \right) \right) \subset \vit \left( t^{-1} \left( t ( P ) \right) \right) = \vit ( P ) .
\end{equation*}
Thus $ \vit \left( t ( P ) \right) \subset t \left( \vit ( P ) \right) $, and so we have equality.
\end{proof}
So the right topology should identify similar point sets so that we aren't really looking at $\ppp(\R^2)$ but rather the quotient of $\ppp(\R^2)$ under similarity.  

On the other hand, forgetting the question of similarity, we do have an intuitive idea of what it means for point sets $P$ and $Q$ in $\ppp(\R^2)$ to be ``close".   In fact there are already metrics, for instance the Hausdorff metric, for this.  The problem is that we need our metric to respect Vorono\"{\i} iteration.  

Suppose we draw little $\epsilon$-balls around all the points of $P$, and discover that each ball contains exactly one point of $Q$.  If so, it is very likely that $\vit(P)$ and $\vit(Q)$ are close as well.   
(The obvious exception is with point sets that have nontrivial cocircularity---jiggling the points a little bit destroys the cocircularity and thus produces extra vertices in the Vorono\"{\i} iteration.)  This idea for a metric is promising but has a serious flaw: it can't compare sets that aren't the same cardinality.   Since there are sets of the same cardinality that iterate to sets of different cardinalities, it would be nice to have the ability to consider whether those iterates are close.

Unfortunately, trying to measure distance between sets of different cardinalities opens Pandora's box.    If several points of $Q$ are inside the $\epsilon$-ball around a point of $P$, we can create all sorts of bizarre patterns in $\vit(Q)$ that may be quite dissimilar to $\vit(P)$.    We find ourselves unsure how to resolve these issues.

\subsection{Are there periodic points?}
\label{ss:perhyp}

For $ k \in \mathbb{ N } $, we say $ P $ is a {\bf period}-$ k $ {\bf point} if $ P \simeq \vit^k ( P ) $.  As a special case, if $ k = 1 $ then $ P \simeq \vit ( P ) $ and we call $ P $ a {\bf fixed point}.  Despite our attempts so far, we have not been able to find any finite periodic point set $ P $.  

\begin{figure}[ht]
\centering
$ \begin{array}{cc}
\includegraphics[width=70mm]{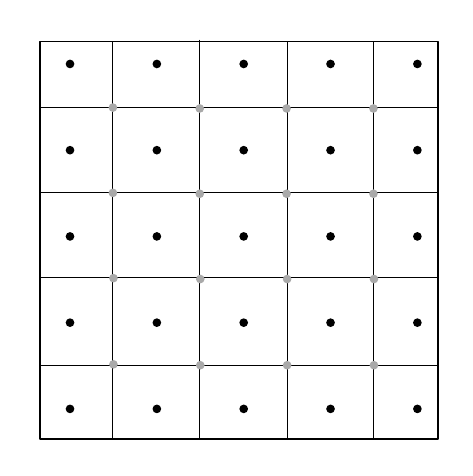}
& 
\includegraphics[width=70mm]{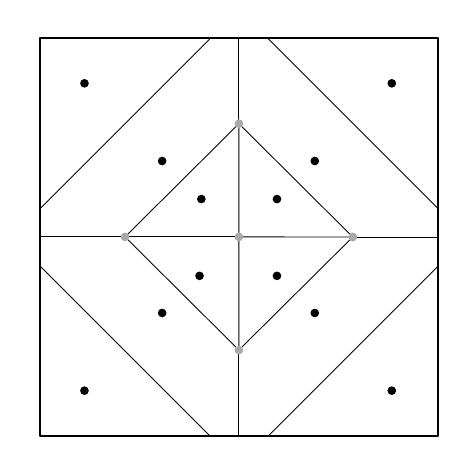} \\
(a) & (b) \\
\end{array} $
\caption{Infinite configurations which are (a) period-$ 1 $ and (b) period-$ 2 $.}
\label{fig:inflattice}
\end{figure}

Curiously, it's not difficult to find an infinite set $ P $ of period $ 1 $ or $ 2 $! In Figure \ref{fig:inflattice} (a), the square lattice iterates to a shifted copy of itself and so has period 1; in Figure \ref{fig:inflattice} (b), points evenly spaced on the diagonals $ y = \pm x $ iterate to points on the $ x $- and $ y $-axes, which then iterate back for a period of $2$.  In both examples one sees a high degree of cocircularity and we assume that this may be a key component in finding examples with other periods as well.

Once a fixed (or periodic) point has been located, we may ask whether it is attracting, repelling, or hyperbolic.  We look at ``nearby'' points (again the need for a well-defined notion of distance!) and look at what their orbits do.  If all nearby points come closer and closer to $ P $, it is attracting; if they all get pushed away we call it repelling.  If there is a mixture it may be hyperbolic.  There's some indication that the grid in Figure \ref{fig:inflattice} (a) may be hyperbolic:  we know it repels grids with ``defects" like the one in Figure \ref{fig:defect}.  However, if we shift an entire row of points up by a fixed amount, the orbit will be pulled back towards the orbit in \ref{fig:inflattice} (a).

\subsection{Sensitive dependence on initial conditions}

\begin{figure}[ht]
\centering
\includegraphics[width=70mm]{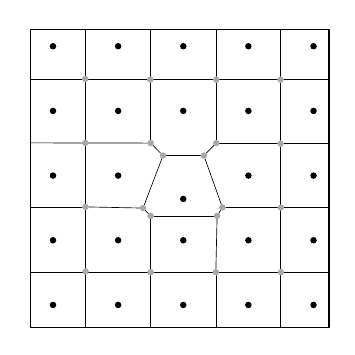}
\caption{Introducing a defect into the square lattice of \ref{fig:inflattice} (a).}
\label{fig:defect}
\end{figure}

Even though we lack a metric on our phase space $ \pha $, there is evidence of sensitive dependence on initial conditions.  For instance, if we shift a point on the grid in Figure \ref{fig:inflattice} (a), it introduces a defect in the Vorono\"{i} iteration, which we see in Figure \ref{fig:defect}.  As the system evolves, that defect will grow to include ever larger regions of the plane until the original grid is no longer recognizable.

One can see this effect in finite configurations as well.  Consider the set depicted in Figure \ref{fig:illuex} (b).  Since $ | \vit(P) | = 4 $ we know it is  doomed to iterate to the empty set.  However, if we were to move any one of the four cocircular points even the slightest amount, we would see the single vertex of degree four become two vertices of degree three.  This five-point configuration has a chance of iterating to five-point sets indefinitely.


\subsection{Can we go backwards?} A sensible question to ask is, what happens when we try to invert the system? Most discrete point sets in the plane are not the vertex set of a Vorono\"{\i} diagram.
So we must ask first when a {\bf preimage} of a point set exists, i.e., a set $ P' $ such that $ \vit ( P' ) = P $.
If there is a preimage, when is it unique?  When it is not unique, how many preimages are possible, and  is there a ``best" one?  

Some literature exists about inverting Vorono\"{\i} tessellations (see e.g. \cite{Ash-Bolker1,inverting}), and we know relatively straightforward conditions that must be satisfied for a preimage to exist.   What is more subtle in our situation is that we don't have a Vorono\"{\i} {\em diagram}---we just have its vertices.   Without the edge adjacencies the problem becomes richer.

It may be that no finite point set exists that can be inverted indefinitely, which would be a shame.   It would be nice if there was a subset of $\ppp(\R^2)$ that was invariant under Vorono\"{\i} iteration in both forward and backwards time.  We conjecture that there is a subset $\sss \subset \ppp(\R^2)$ of five-point configurations with 2 in the interior and 3 on the boundary for which (1) the set is invariant under forward Vorono\"{\i} iteration, and (2) there exist preimages of all orders within the set.    If we could identify criteria for membership in $\sss$, then we could begin to examine the orbits of $\sss$ for evidence of chaos.

\subsection{Conclusion}
We admit that we don't know much about the behavior of this dynamical system yet.  We know exactly what happens for small point sets, but as soon as there are five or more points we encounter difficulties.  We have a theorem \cite{efgh} that measures how the point sets grow or shrink in size, but it requires geometric information that isn't always readily available.  We have evidence that many point sets grow without bound, but have been unable to determine which conditions guarantee this.   We've noticed evidence of sensitive dependence on initial conditions, but we don't have the machinery to measure the phenomenon.  In the course of our study we've developed tools to help us with our insight on this finicky dynamical system, and we've shared some of that here.   There is one thing we know for certain: plenty of problems remain.  Some are simple enough to be studied by budding mathematicians, and some may be subtle enough to interest their advisers as well. 

\bibliographystyle{amsplain}

\end{document}